\documentclass[journal,12pt,onecolumn]{IEEEtran}

\ifCLASSINFOpdf

\else

\fi

\usepackage{graphics} 
\usepackage{epsfig} 
\usepackage{mathptmx} 
\usepackage{times} 
\usepackage{amsmath} 
\usepackage{amssymb} 
\usepackage{graphicx,paralist}
\usepackage{float}
\usepackage[T1]{fontenc}
\usepackage{array}
\usepackage{color}
\usepackage{xfrac}
\usepackage{cite}
\usepackage{textcomp}
\usepackage{empheq}
\usepackage{tikz}
\usepackage{enumitem}

\begin{document}

\title{Subsystem-Based Control with Modularity for Strict-Feedback Form Nonlinear Systems}
\author{Janne Koivum\"{a}ki$^1$, Jukka-Pekka Humaloja$^2$, Lassi Paunonen$^2$, Wen-Hong 
Zhu$^3$\\ and Jouni Mattila$^1$
\thanks{$^1$J. Koivum\"{a}ki and J. Mattila are with Automation and Mechanical Engineering, Faculty 
of Engineering and Natural Sciences, Tampere University, Finland. e-mail: firstname.surname@tuni.fi}
\thanks{$^2$J.-P. Humaloja and L. Paunonen are with Mathematics, Faculty of Information Technology and Communication Sciences, Tampere University, Finland. e-mail: firstname.surname@tuni.fi}
\thanks{$^3$W.-H. Zhu is with the Canadian Space Agency, Canada. e-mail: WenHong.Zhu@canada.ca}
}

\maketitle 

\begin{abstract}
This study proposes an adaptive \textit{subsystem-based control} (SBC) for systematic and straightforward nonlinear~control of \textit{n}th-order strict-feedback form (SFF) systems.~By decomposing the SFF system to subsystems, a generic~term~(namely \textit{stability connector}) can be created to address dynamic interactions between the subsystems. This 1) enables modular control design with global asymptotic stability, 2) such that both the control design and the stability analysis can be performed locally at a subsystem level, 3) while avoiding an excessive growth of the control design complexity when the system order \textit{n} increases. The latter property makes the method suitable especially for high-dimensional systems. We also design a smooth projection function for addressing system parametric uncertainties. Numerical simulations demonstrate the efficiency of the method.
\end{abstract}


\IEEEpeerreviewmaketitle

\section{Introduction}

\IEEEPARstart{N}{onlinear} model-based control aims to design a specific feedforward (FF) 
compensation term based on the system inverse dynamics to generate the control output(s)~from 
the system states and desired input signals \cite{Slotine1991}. If the FF compensation can exactly 
capture the inverse of the plant dynamics for all frequencies, an infinite control bandwidth~with zero 
tracking error becomes theoretically possible \cite{Skogestad2007,Zhu_VDC}. While early control 
methods, e.g., \textit{feedback linearization} \cite{Hunt1983}, aimed to cancel (or linearize) the 
system nonlinearities, \textit{adaptive backstepping} \cite{Backstepping1995} became a significant 
breakthrough in nonlinear systems control by incorporating the nonlinearities towards ideal FF 
compensation with \textit{global asymptotic stability}. 

This study proposes globally asymptotically stable adaptive \textit{subsystem-based control} (SBC) 
for \textit{n}th-order strict-feedback form (SFF) systems. The proposed method has built-in 
modularity and it avoids excessive growth of the control design complexity when the system order 
\textit{n} increases (an issue reported for backstepping-based methods in several studies 
\cite{Swaroop2000_DSC,Song2011_DSC,Yip1998,Wang2005}). \textit{Dynamic surface control} 
(DSC) \cite{Swaroop2000_DSC,Song2011_DSC} and \textit{adaptive} DSC \cite{Yip1998} are 
previously developed as an alternative to backstepping to avoid the reported ``explosion of 
complexity'' with semi-global stability. They are based on \textit{multiple sliding surface} (MSS) 
control \cite{Slotine1993,Won1996} (a method similar~to~backstepping) using a series of low-pass 
filters \cite{Song2011_DSC}. Our method does not employ filtering and achieves global asymptotic 
stability.

The proposed method originates from \textit{virtual decomposition control} (VDC) \cite{Zhu_VDC}, 
\cite{Zhu1997} that is developed~for controlling complex robotic systems. \textit{Modularity} is 
one of the key aspects in addressing complexity in advanced control 
realizations \cite{Mastellone2021}, \cite[Sec. IV]{Mattila_TMECH2017}. In VDC, robotic systems are 
virtually decomposed into \textit{modular subsystems} (rigid links and joints) such that both control 
design and stability analysis~can~be~performed locally at the subsystem (SS) level to 
guarantee overall global asymptotic stability. In particular, VDC introduced \textit{virtual 
power flows} (VPFs) \cite[Def. 2.16]{Zhu_VDC} to define dynamic interactions between the adjacent 
SSs such that the VPFs cancel each others out when the SSs are connected. However,~when~applied 
beyond robotics, the interactions between SSs will no longer be described by VPFs 
\cite{Koivumaki_ASME2017}. Some~early~ideas for the proposed method originate from the 
application-oriented paper in \cite{Koivumaki_ASME2017}. In addition, some ideological similarities can be seen to the passivity-based approach in \cite{SeiAll02} for controlling SFF systems with global asymptotic stability. While the method in \cite{SeiAll02} designed strictly passive interaction dynamics for adjacent SSs, we propose new generic tools to compensate the interaction dynamics such that every SS is automatically stabilized by its adjacent SS. More details on differences to \cite{SeiAll02} can be found in Remark~\ref{diff_to_SeyAll02}.

As the main contribution, the proposed method generalizes the ``subsystem-based control philosophy'' in \cite{Zhu_VDC,Koivumaki_ASME2017}~for~controlling the \textit{n}th-order SFF systems. After defining a generic form for SSs, we design a specific \textit{stability connector} (a generic spill-over term in SS stability analysis in Def. \ref{def:stab_con})~to address dynamic interactions between the adjacent~SSs. We show that every SS with a ``stability preventing'' connector is compensated by the subsequent SS with a corresponding ``stabilizing'' connector. Similarly to VDC, we formulate a generic definition for \textit{virtual stability}\footnote{In terms of Lyapunov functions, definition of \textit{virtual stability} (see Def. \ref{def:Vstab} in Section \ref{sec:stability}) includes quadratic terms for asymptotic convergence added with stability connector(s) for compensating/stabilizing dynamics of adjacent SSs.} such that when every SS is virtually stable, the overall system becomes automatically globally asymptotically stable. Instead of using~Lebesque~$L_2$/$L_\infty$ integrable functions as in \cite{Zhu_VDC, SeiAll02, Koivumaki_ASME2017}, we base the results on Lyapunov functions. The proposed method is modular in the sense that control laws for every SS can be designed with a single generic-form equation as shown in Remarks~\ref{remark:modularity} and \ref{remark:modularity_adaptive}. As part of the control design, we design a smooth projection function to address the system parametric uncertainties.

Next, Section \ref{sec:problem_statement} introduces the control problem.~Section~\ref{sec:ASDB_control} formulates the proposed method. Section \ref{sec:stability} provides in-depth analysis on the control design and its stability.~Section~\ref{sec:validation}~provides numerical validation. Section \ref{sec:conclusions} concludes the study.  

\newtheorem{theorem}{Theorem}[section]
\newtheorem{lemma}{Lemma}[section]
\newtheorem{proposition}{Proposition}[section]
\newtheorem{assump}{Assumption}[section]
\newtheorem{defin}{Definition}[section]
\newtheorem{remark}{Remark}[section]

\vspace{-0.00cm}
\section{The Control Problem}
\label{sec:problem_statement}

Consider the following \textit{n}th-order SFF system
\begin{empheq}[left=\empheqlbrace]{align}
  \theta_{11}\dot{x}_1 &= f_1({x}_1) + g_1({x}_1)x_2 \label{EQ_dyn1}\\
  \theta_{i1}\dot{x}_i &= f_i(\pmb{x}_i) + g_i(\pmb{x}_i)x_{i+1},\ \forall i \in \{2,...,n-1\} 	\label{EQ_dyn2}\\
  \theta_{n1}\dot{x}_n &= f_n(\pmb{x}_n) + g_n(\pmb{x}_n)u \label{EQ_dyn3} 
\end{empheq}
where $\pmb{x}_k = [x_1,x_2,\cdots,x_k]$ for all $k \in \{1,\cdots,n\}$, $u$ is the system input, $f_k(\pmb{x}_k)$ for all $k \in \{1,\cdots,n\}$ can be further written as
\begin{equation}
	f_k(\pmb{x}_k) = \theta_{k2}\gamma_{k2}(\pmb{x}_k) + \theta_{k3}\gamma_{k3}(\pmb{x}_k) + \cdots + \theta_{kj}\gamma_{kj}(\pmb{x}_k) \label{EQ_lin_par}
\end{equation}
and $\theta_{k1},\theta_{k2},\cdots,\theta_{kj} > 0$ in \eqref{EQ_dyn1}--\eqref{EQ_lin_par} are the 
system parameters. Similarly to backstepping, we assume that $g_k(t,\pmb{x}_k)$~and 
$f_k(t,\pmb{x}_k)$ (i.e, $\gamma_{k\zeta}(t,\pmb{x}_k)$, $\forall \zeta \in \{2,\ldots,j\}$) are 
sufficiently smooth and $g_k(t,\pmb{x}_k) \neq 0$ on $[0, \infty) \times \mathbb{R}^k$.

Throughout the paper, we use \textit{n} to denote the system overall order, while it also denotes the 
last SS (or its element) in \eqref{EQ_dyn3}. We use $i \in \{2,\cdots,n-1\}$ to denote a SS (or its 
element)~in~the middle of the SFF sequence; see \eqref{EQ_dyn2}. We use $k$ to denote an 
arbitrary decomposed SS (or its element), such that generic form for the \textit{k}th SS (i.e., 
SS$_k$) in \eqref{EQ_dyn1}--\eqref{EQ_dyn3} is given by
%An arbitrary decomposed SS (i.e., SS$_k$) in \eqref{EQ_dyn1}--\eqref{EQ_dyn3} can be defined~as
\begin{equation}
\theta_{k1}\dot{x}_k = f_k(\pmb{x}_k) + g_k(\pmb{x}_k)x_{k+1},\ \forall k \in \{1,...,n\} \label{EQ_generic_SS}
\end{equation}
where we denote $x_{n+1} = u$.

Let $x_{{1}\rm d}(t) \in C^{n-1}(0,\infty)$ be a desired trajectory for $x_{{1}}(t)$ such that $x_{1{\rm d}}^{(n)}$ exists almost everywhere. Next, our aim is to design~a control for the system in \eqref{EQ_dyn1}--\eqref{EQ_dyn3}, such that $e_1(t)$ = $x_{{1}\rm d}(t) - x_1(t)$ globally asymptotically converges to zero when $t>0$.

\section{The Proposed Control Method}
\label{sec:ASDB_control}

In Section \ref{subsec:SBC_design}, we first design the baseline~SBC~by~assuming the plant parameters ${\theta}_{kj}$ in \eqref{EQ_dyn1}--\eqref{EQ_lin_par} known~$\forall k$,~$\forall j$.~Then, Section \ref{subsec:projF} proposes a projection function $\mathcal{P}_{k}$~for parametric uncertainties, such that SBC can be updated to the proposed \textit{adaptive} SBC in Section \ref{subsec:ASBC_design}. The control design philosophy behind the proposed method is analyzed later~in~Section \ref{sec:stability}.

\vspace{-0.0cm}
\subsection{Subsystem-Based Control}
\label{subsec:SBC_design}

Assume that the system in \eqref{EQ_dyn1}--\eqref{EQ_lin_par} is not subject to~any~parametric uncertainty in ${\theta}_{kj}$, $\forall k, \forall j$. The baseline SBC for the SFF system in \eqref{EQ_dyn1}--\eqref{EQ_dyn3} can be designed as
\begin{empheq}[left=\empheqlbrace]{align}
	g_1({x}_1)x_{\rm 2d}	&= {\theta}_{11}\dot{x}_{1{\rm d}} - {f}_1({x}_1) + \lambda_1{e}_1 \nonumber\\ 
												&={\mathbf Y}_{1}{\pmb \theta}_1 + \lambda_1{e}_1 \label{ctrl1}\\
	g_i(\pmb{x}_i)x_{(i+1){\rm d}} &= {\theta}_{i1}\dot{x}_{i{\rm d}} - {f}_i(\pmb{x}_i) + \delta_{i-1}g_{i-1}(\pmb{x}_{i-1}){e}_{i-1} + \lambda_i{e}_{i} \nonumber\\
												&= {\mathbf Y}_{i}{\pmb \theta}_i + \delta_{i-1}g_{i-1}(\pmb{x}_{i-1}){e}_{i-1} + \lambda_i{e}_{i} \label{ctrl2}\\
	g_n(\pmb{x}_n)u &= {\theta}_{n1}\dot{x}_{n{\rm d}} - {f}_n(\pmb{x}_n) + \delta_{n-1}g_{n-1}(\pmb{x}_{n-1}){e}_{n-1} + \lambda_n{e}_{n} \nonumber\\
												&= {\mathbf Y}_{n}{\pmb \theta}_n + \delta_{n-1}g_{n-1}(\pmb{x}_{n-1}){e}_{n-1} + \lambda_n{e}_{n} \label{ctrl3}
\end{empheq}
where $\lambda_k e_k = \lambda_k({x}_{k{\rm d}} - {x}_{k})$ is the \textit{local feedback} (FB) \textit{term} with $\lambda_k > 0$; $\delta_{k-1}g_{k-1}(\pmb{x}_{k-1}){e}_{k-1}$ is the \textit{stabilizing FB term} for the previous subsystem, $\delta_{k-1} > 0$; ${f}_k(\pmb{x}_k)$ is defined in \eqref{EQ_lin_par}; and in the model-based FF compensation term ${\mathbf Y}_{k}{\pmb \theta}_k$, the regressor ${\mathbf Y}_{k}$ and the parameter vector ${\pmb \theta}_k$ are defined~as
\begin{align}
	{\mathbf Y}_{k} & := \big[\dot{x}_{k{\rm d}},\ -\gamma_{k2}(\pmb{x}_k),\ 
	-\gamma_{k3}(\pmb{x}_k), \cdots,\ -\gamma_{kj}(\pmb{x}_k)\big] \in \mathbb{R}^{1\times j} 
	\label{EQ_Yk}\\ 
	{\pmb \theta}_k & := \left[\theta_{k1},\theta_{k2},\theta_{k3},\cdots,\theta_{kj}\right]^T \in 
	\mathbb{R}^{j}. \label{EQ_thetak}
\end{align}   

Similarly to backstepping, $x_{(k+1){\rm d}}$ in \eqref{ctrl1} and \eqref{ctrl2} acts~as~a fictitious control from SS$_{k}$ to the subsequent SS, $\forall k \in \{1,...,n-1\}$. The real control effort $u$ can be obtained from \eqref{ctrl3} after stepping through every SS.
\begin{remark} \label{remark:modularity}
Similarly to SS$_k$ dynamics in \eqref{EQ_generic_SS}, the control~in \eqref{ctrl1}--\eqref{ctrl3} can be reproduced with a generic and modular equation  
\begin{align} \nonumber%\label{EQ_SSk_ctrl}
		g_k(\pmb{x}_k)x_{(k+1){\rm d}} = {\mathbf Y}_{k}{\pmb \theta}_k + \delta_{k-1}g_{k-1}(\pmb{x}_{k-1}){e}_{k-1} + \lambda_k{e}_{k}
\end{align}
 $\forall k \in \{1,...,n\}$, such that $\delta_{0}g_{0}(\pmb{x}_{0}){e}_{0} = 0$ and $x_{(n+1){\rm d}} = u$. The modularity in the control provides~that~changing SS$_k$ dynamics, or adding/removing SSs, do not alter the structure of control laws in the remaining~SSs.
\end{remark}

\vspace{-0.0cm}
\subsection{The Proposed Smooth Projection Function}
\label{subsec:projF}

\begin{defin} \label{def:PA3}
A piecewise-continuous function $\mathcal{P}_{k}({\rm p}(t),\rho,$ $\sigma,a,b,c,t) \in \mathbb{R}$ is a \textit{k}th-order differentiable scalar function, $\forall k \in \{1,...,n\}$, defined for $t \geqslant 0$ such that its time derivative is governed by
\begin{equation} \label{EQ_PAdef3}
\dot{\mathcal{P}}_{k} = \rho\left({\rm p}(t) + \sigma{\kappa}\right)
\end{equation}
where $\rho,\sigma>0$, ${\rm p}(t) \in C^{n-k}(0,\infty; \mathbb{R})$, $\forall k \in \{1,...,n\}$, and
\begin{align} \nonumber%\label{EQ_kappa3}
{\kappa} = \left\{
  \begin{array}{l l l}
	  (b - \mathcal{P}_{k}), & \text{if}\hspace{1.36cm} \mathcal{P}_{k} \geqslant b + c\\
		(b - \mathcal{P}_{k})S_b(\mathcal{P}_{k}), & \text{if}\hspace{0.75cm} b < \mathcal{P}_{k} < b + c\\
		0, & \text{if}\hspace{0.75cm} a \leqslant \mathcal{P}_{k} \leqslant b\\
		(a - \mathcal{P}_{k})S_a(\mathcal{P}_{k}), & \text{if}\hspace{0.22cm} a - c < \mathcal{P}_{k} < a\\
    (a - \mathcal{P}_{k}), & \text{if}\hspace{1.36cm} \mathcal{P}_{k} \leqslant a - c
		\end{array} \right.
\end{align}
where $a,b,c > 0$ satisfy $c + b > b \geqslant a > a-c > 0$; $S_a(\mathcal{P}_{k}) \in C^{n-k}: 
(a-c, a) \to (1, 0)$ is \textit{strictly decreasing}; and $S_b(\mathcal{P}_{k}) \in C^{n-k}: (b, b+c) \to 
(0, 1)$ is \textit{strictly increasing}.
\end{defin}

A solution for the switching functions $S_a(\mathcal{P}_{k})$ and $S_b(\mathcal{P}_{k})$ can be found in Appendix \ref{App:projF} that also provides a detailed~analysis on the projection function $\mathcal{P}_{k}$ and its properties.

\vspace{-0.0cm}
\subsection{Adaptive Subsystem-Based Control}
\label{subsec:ASBC_design}

Let the system in \eqref{EQ_dyn1}--\eqref{EQ_lin_par} be subject to parametric uncertainties, i.e., ${\theta}_{kj}$ is unknown $\forall k, \forall j$. The control in Section~\ref{subsec:SBC_design} can be updated to the proposed \textit{adaptive} SBC as 
\begin{empheq}[left=\empheqlbrace]{align}
	g_1({x}_1)x_{\rm 2d}	&={\mathbf Y}_{1}\widehat{\pmb \theta}_1 + \lambda_1{e}_1 \label{EQ_ctrl1a}\\
	g_i(\pmb{x}_i)x_{(i+1){\rm d}} &= {\mathbf Y}_{i}\widehat{\pmb \theta}_i + \delta_{i-1}g_{i-1}(\pmb{x}_{i-1}){e}_{i-1} + \lambda_i{e}_{i} \label{EQ_ctrl2a}\\
	g_n(\pmb{x}_n)u &= {\mathbf Y}_{n}\widehat{\pmb \theta}_n + \delta_{n-1}g_{n-1}(\pmb{x}_{n-1}){e}_{n-1} + \lambda_n{e}_{n} \label{EQ_ctrl3a}
\end{empheq}
where ${\mathbf Y}_{k}\widehat{\pmb \theta}_k$ is the adaptive model-based~FF compensation, ${\mathbf Y}_{k}$ is defined in \eqref{EQ_Yk} and $\widehat{\pmb \theta}_k \in \mathbb{R}^{j}$ is an estimate of ${\pmb \theta}_k$ in \eqref{EQ_thetak}.~The estimated parameters in $\widehat{{\pmb \theta}}_k$ need to be updated. We define
\begin{align}
	{\mathbf{p}}_{k} := {e}_k{{\mathbf Y}}^T_{k} \label{PA_pk} 
\end{align}
such that the $\zeta$th element of $\widehat{{\pmb\theta}}_{k}$ in \eqref{EQ_ctrl1a}--\eqref{EQ_ctrl3a} can be updated by using the projection function $\mathcal{P}_k$ in Definition \ref{def:PA3} as
\begin{align}
	\widehat{\theta}_{k\zeta} = \mathcal{P}_k({\rm p}_{k\zeta},\rho_{k\zeta},\sigma_{k\zeta},\underline{\theta}_{k\zeta},\overline{\theta}_{k\zeta},c_{k\zeta},t), \forall \zeta \in \{1,...,j\} \label{PA_thetak}
\end{align}
where $\widehat{\theta}_{k\zeta}$ is the $\zeta$th element of $\widehat{{\pmb\theta}}_{k}$; ${\rm p}_{k\zeta}$ is the $\zeta$th element of ${\mathbf{p}}_{k}$ in \eqref{PA_pk}; $\rho_{k\zeta} > 0$ and $\sigma_{k\zeta} > 0$ are the parameter update gains; $\underline{\theta}_{k\zeta}$ and $\overline{\theta}_{k\zeta}$ are the lower and the upper bounds of ${\theta}_{k\zeta}$; and $c_{k\zeta}$ defines the activation interval beyond the bounds. 

\begin{figure}
	\vspace{0.0cm}
	\centering
	\includegraphics[width=1.00\textwidth]{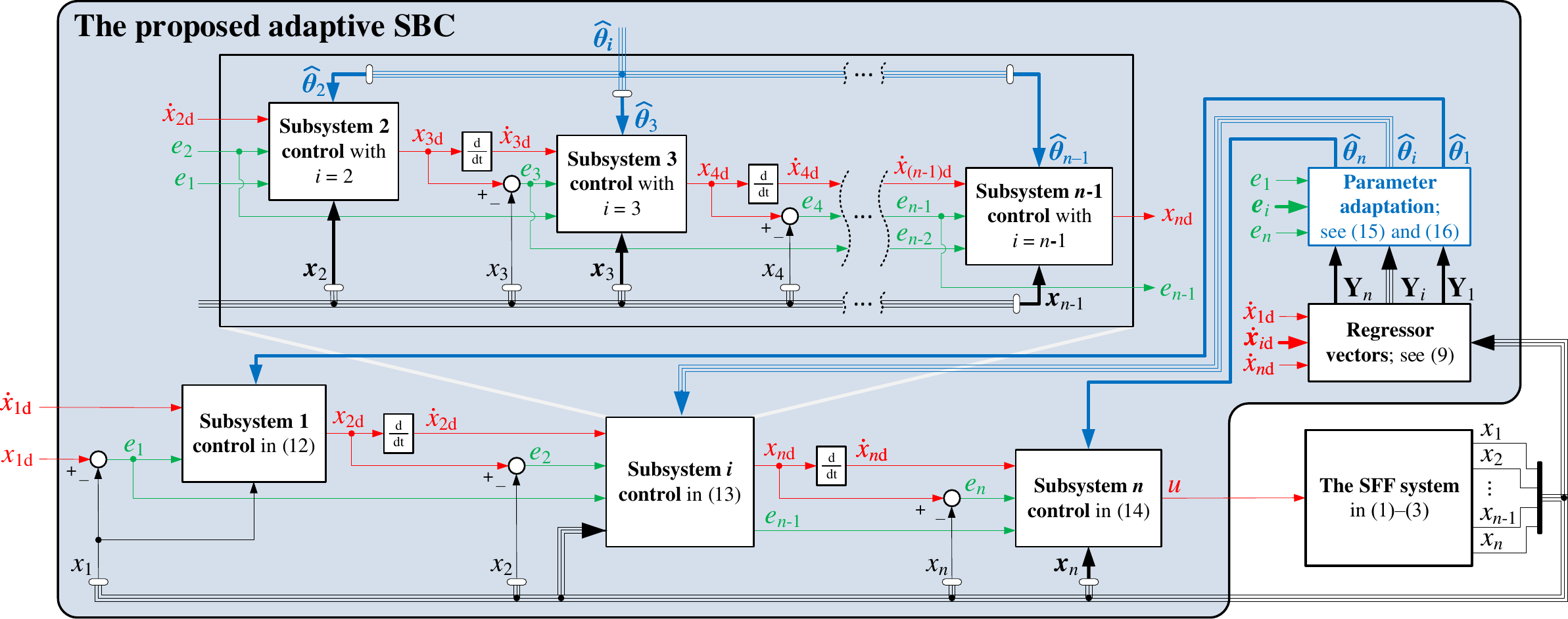}
	\vspace{-0.7cm}
	\caption{Diagram of the proposed adaptive SBC (highlighted in light blue). The desired variables (and the control output $u$) are shown in red, the feedback~signals are in green, the adaptive control is in blue, and the system output states are in black. The bold lines are vectors and the thin lines are scalar variables.}
	\label{fig:ctrl_design}
	\vspace{-0.0cm}
\end{figure}

Fig. \ref{fig:ctrl_design} shows the diagram of the proposed method. 

\begin{remark}
As Fig. \ref{fig:ctrl_design} and \eqref{EQ_ctrl1a}--\eqref{EQ_ctrl3a} show, $\widehat{\pmb \theta}_k$~in~SS$_k$~should be continuously differentiable in $C^{n-k}$ when stepping through the remaining SSs. The projection function $\mathcal{P}_k$ in \eqref{EQ_PAdef3} satisfies $\widehat{\pmb \theta}_k \in C^{n-k}$, $\forall k \in \{1,...,n\}$. If SS$_k$ satisfies $n-k \leqslant 1$,~a~projection function $\mathcal{P} \in C^1$ in \cite{Zhu_VDC,Zhu1999adaptive} can be used instead of $\mathcal{P}_k$. If SS$_k$ satisfies $n-k \leqslant 2$, $\mathcal{P}_2 \in C^2$ in \cite{Zhu_VDC,Zhu1999adaptive} can be used.
\end{remark}
\begin{remark} \label{remark:modularity_adaptive}
As in Remark \ref{remark:modularity}, SS$_k$ control in \eqref{EQ_ctrl1a}--\eqref{EQ_ctrl3a}~can be reproduced with a generic and modular equation    
\begin{align} \nonumber%\label{EQ_SSk_ctrla}
		g_k(\pmb{x}_k)x_{(k+1){\rm d}} = {\mathbf Y}_{k}\widehat{\pmb \theta}_k + \delta_{k-1}g_{k-1}(\pmb{x}_{k-1}){e}_{k-1} + \lambda_k{e}_{k}
\end{align}
$\forall k \in \{1,...,n\}$, such that $\delta_{0}g_{0}(\pmb{x}_{0}){e}_{0} = 0$ and $x_{(n+1){\rm d}} = u$. The modularity in the control~provides~that~changing SS$_k$ dynamics, or adding/removing SSs, do not alter the structure of control laws in the remaining~SSs.
\end{remark}
\begin{remark} \label{diff_to_SeyAll02}
As the main difference to \cite{SeiAll02},~we~design~stabilizing FB term $\delta_{k-1}g_{k-1}(\pmb{x}_{k-1}){e}_{k-1}$, $\forall k \in \{2,...,n\}$, in \eqref{EQ_ctrl1a}--\eqref{EQ_ctrl3a} to produce \textit{stability connector} $s_{k-1}$ (analyzed next~in~Section \ref{sec:stability}), such that passivity between SSs do not need to be considered. While the results in \cite{SeiAll02} are based on Lebesque $L_2$/$L_\infty$ integrable functions, we base the results on Lyapunov functions. We also proposed novel projection function $\mathcal{P}_k$ in Definition \ref{def:PA3} to address the system parametric uncertainties.  
\end{remark}

\vspace{-0.0cm}
\section{Stability Analysis}
\label{sec:stability}
 \vspace{-0.0cm}

Next, we provide an in-depth analysis on the adaptive SBC in Section \ref{subsec:ASBC_design}. Respective analysis can be performed for the SBC in Section \ref{subsec:SBC_design} using $\pmb{\theta}_k - {\pmb{\theta}}_k = 0$ instead of $\pmb{\theta}_k - \widehat{\pmb{\theta}}_k$.

Motivated by a key concept in virtual stability analysis---a \textit{virtual power flow}
\cite[Sect. 2.9.2]{Zhu_VDC}---we introduce a related notion of a \textit{stability connector} as follows:
\begin{defin} \label{def:stab_con}
For the system \eqref{EQ_dyn1}--\eqref{EQ_dyn3} with the control \eqref{EQ_ctrl1a}--\eqref{EQ_ctrl3a}, the \textit{stability connector} $s_{k}$ is defined~as
\begin{align}
	s_{k} = \Delta_kg_k(t,\pmb{x}_k){e}_{k}{e}_{k+1} \nonumber
\end{align}	
where SS-related term $\Delta_k = 1$, if $k=1$, and $\Delta_k = \frac{1}{\delta_1\cdots\delta_{k-1}}$, if $k > 1$, and $\delta_1,\delta_2,\cdots,\delta_{k-1}>0$ are feedback gains from Section~\ref{sec:ASDB_control}.
\end{defin}

Next, in Lemmas \ref{lem:ss_1a}--\ref{lem:ss_na} we provide auxiliary results for the convergence analysis in Theorem \ref{thm:stability}. Motivated by the concept of \textit{virtual stability} \cite[Sect. 2.9]{Zhu_VDC}, the auxiliary analysis is carried out for the individual subsystem error dynamics $e_k$ and the corresponding parameter estimation errors $\pmb{\theta}_k - \widehat{\pmb{\theta}}_k$. 

Subtracting \eqref{EQ_dyn1} from \eqref{EQ_ctrl1a}, adding $\theta_{11}\dot{x}_{\rm 1d} - \theta_{11}\dot{x}_{\rm 1d}$ = 0, using \eqref{EQ_lin_par}, \eqref{EQ_Yk} and \eqref{EQ_thetak}, and rearranging the terms, we get~the~following error dynamics for SS$_1$ 
\begin{align}
	\theta_{11}\dot{e}_1 = - \lambda_1{e}_1 + g_1({x}_1)e_{\rm 2} + {\mathbf Y}_{1}({\pmb \theta}_1 - \widehat{\pmb \theta}_1). \label{EQ_error_dyn1a}
\end{align}

\begin{lemma} \label{lem:ss_1a} 
Considering SS$_1$ error dynamics in \eqref{EQ_error_dyn1a}, and $\pmb{\theta}_1 - \widehat{\pmb{\theta}}_1$~governed by \eqref{EQ_thetak}, \eqref{PA_pk} and \eqref{PA_thetak},~the derivative of the quadratic function
\begin{align}
		\nu_{1} = \frac{1}{2}\left(\theta_{11}{e}_1^2 + \sum_{\zeta=1}^{j}\frac{(\theta_{1\zeta} - \widehat{\theta}_{1\zeta})^2}{\rho_{1\zeta}}\right) \label{EQ_nu1a}
\end{align} 
along the trajectories of the error dynamics satisfies
\begin{align}
	\dot{\nu}_{1} \leqslant - \lambda_1{e}_1^2 + s_1 \label{EQ_nu1a_dot}
\end{align} where $s_1$ is the stability connector from Definition \ref{def:stab_con}.
\end{lemma}

\begin{IEEEproof}
See Appendix \ref{App:lemmas_ss}.
\end{IEEEproof}

\begin{remark}
In Lemma \ref{lem:ss_1a}, term $e_2$ in \eqref{EQ_error_dyn1a} is treated as an external input that causes $s_1$ to appear in \eqref{EQ_nu1a_dot} (see Appendix~\ref{App:lemmas_ss}) that will be canceled out based on the result of the next lemma. The dynamics of $e_2$ as well as the subsequent subsystems error dynamics are accounted for in the next~two~lemmas. 
\end{remark}

Subtracting \eqref{EQ_dyn2} from \eqref{EQ_ctrl2a}, adding $\theta_{i1}\dot{x}_{i{\rm d}} - \theta_{i1}\dot{x}_{i{\rm d}}$ = 0 using \eqref{EQ_lin_par}, \eqref{EQ_Yk} and \eqref{EQ_thetak}, and rearranging the terms, we get~the~following error dynamics for SS$_i$, $\forall i \in \{2,...,n-1\}$,
\begin{align}
  \theta_{i1}\dot{e}_i &= - \lambda_i{e}_{i} - \delta_{i-1}g_{i-1}(\pmb{x}_{i-1}){e}_{i-1} + g_i(\pmb{x}_i)e_{i+1} \nonumber\\
  &\hspace{0.4cm}+{\mathbf Y}_i({\pmb \theta}_i - \widehat{\pmb \theta}_i). \label{EQ_error_dyn2a}
\end{align}

\begin{lemma} \label{lem:ss_ia}
Considering SS$_i$ error dynamics in \eqref{EQ_error_dyn2a}, and $\pmb{\theta}_i - 
\widehat{\pmb{\theta}}_i$ governed by \eqref{EQ_thetak}, \eqref{PA_pk} and \eqref{PA_thetak}, the 
derivative of the quadratic function
\begin{equation}
	\nu_{i} =
	\frac{1}{2(\delta_1\cdots\delta_{i-1})}\left(\theta_{i1}{e}_{i}^2
	+
	\sum_{\zeta=1}^{j}\frac{(\theta_{i\zeta}
	- \widehat{\theta}_{i\zeta})^2}{\rho_{i\zeta}} \right) \label{EQ_nu2a}
\end{equation} along the trajectories of the error dynamics satisfies
\begin{align}
	\dot{\nu}_{i} \leqslant -
	\frac{\lambda_i}{\delta_1\cdots\delta_{i-1}}{e}_{i}^2
	- s_{i-1} + s_{i} \label{EQ_nu2a_dot}
\end{align} where $s_{i-1}$ and $s_i$ are the stability connectors from Definition \ref{def:stab_con}.
\end{lemma}

\begin{IEEEproof}
%See Appendix \ref{App:lem:ss_ia}.
See Appendix \ref{App:lemmas_ss}.
\end{IEEEproof}

\begin{remark}
Similarly to Lemma \ref{lem:ss_1a}, $e_{i+1}$ in \eqref{EQ_error_dyn2a} is treated as an external input that causes $s_i$ to appear in \eqref{EQ_nu2a_dot}. The stabilizing FB term $\delta_{i-1}g_{i-1}(\pmb{x}_{i-1}){e}_{i-1}$ in \eqref{EQ_error_dyn2a} creates another stability connector $-s_{i-1}$ to appear in \eqref{EQ_nu2a_dot} (see Appendix~\ref{App:lemmas_ss}) that will cancel out $s_{i-1}$ from the previous SS. The last connector $s_{n-1}$ will be canceled out based on the result of the next lemma, after which we are in the position to present the convergence result for the overall error dynamics.
\end{remark}

Subtracting \eqref{EQ_dyn3} from \eqref{EQ_ctrl3a}, adding $\theta_{n1}\dot{x}_{n{\rm d}} - \theta_{n1}\dot{x}_{n{\rm d}}$ = 0 using \eqref{EQ_lin_par}, \eqref{EQ_Yk} and \eqref{EQ_thetak}, and rearranging the terms, we get~the~following error dynamics for SS$_n$ 
\begin{align}
  \theta_{n1}\dot{e}_n = - \lambda_n{e}_{n} - \delta_{n-1}g_{n-1}(\pmb{x}_{n-1}){e}_{n-1} + {\mathbf Y}_{n}({\pmb \theta}_n - \widehat{\pmb \theta}_n). \label{EQ_error_dyn3a}
\end{align}

\begin{lemma} \label{lem:ss_na}
Considering SS$_n$ error dynamics in \eqref{EQ_error_dyn3a}, and $\pmb{\theta}_n - \widehat{\pmb{\theta}}_n$ governed by \eqref{EQ_thetak}, \eqref{PA_pk} and \eqref{PA_thetak},~the derivative of the quadratic function 
\begin{equation}
	\nu_{n} =
	\frac{1}{2(\delta_1\cdots\delta_{n-1})}\left( \theta_{n1}{e}_{n}^2
	+ \sum_{\zeta=1}^{j}\frac{(\theta_{n\zeta}
	- \widehat{\theta}_{n\zeta})^2}{\rho_{n\zeta}}  \right)\label{EQ_nu3a}
    \end{equation}
    along the trajectories of the error dynamics satisfies
\begin{align}
	\dot{\nu}_{n} \leqslant -
	\frac{\lambda_n}{\delta_1\cdots\delta_{n-1}}{e}_{n}^2
	- s_{n-1} \label{EQ_nu3a_dot}
\end{align} 
where $s_{n-1}$ is the stability connector from Definition \ref{def:stab_con}.
\end{lemma}

\begin{IEEEproof}
%See Appendix \ref{App:lem:ss_na}.
See Appendix \ref{App:lemmas_ss}.
\end{IEEEproof}

We will now construct a Lyapunov candidate for the overall error dynamics as the sum of the quadratic functions from Lemmas \ref{lem:ss_1a}--\ref{lem:ss_na}. Based on the properties derived in the lemmas, we obtain that the error dynamics will remain bounded, and moreover, that the control errors converge globally asymptotically to zero. The result is given in the following theorem.
\begin{theorem} \label{thm:stability}
  Consider the error dynamics $\pmb{e} = [e_1, \ldots, e_n]^T$ and the parameter estimation error $\pmb{\theta}_k - \widehat{\pmb{\theta}}_k$, $\forall k \in \{1,2,\ldots, n\}$, that are governed in Lemmas \ref{lem:ss_1a}--\ref{lem:ss_na}. For arbitrary initial conditions, $\pmb{\theta}_k - \widehat{\pmb{\theta}}_k$ remains bounded and $e_k(t) \to 0$ globally as $t \to \infty$ for all $k \in \{1,2,\ldots, n\}$.
\begin{IEEEproof}
Using \eqref{EQ_nu1a}, \eqref{EQ_nu2a} and \eqref{EQ_nu3a}, we choose a Lyapunov candidate function for the overall error dynamics as
\begin{align} 
\nu_{\rm tot} &= \nu_1 + \sum_{i = 2}^{n
- 1} {\nu_i} + \nu_n \nonumber\\
&= \frac{1}{2}\pmb{e}^T{\mathbf A}\pmb{e} + \sum_{k = 1}^n\frac{1}{2(\delta_1\cdots\delta_{k-1})}\sum_{\zeta=1}^{j}\frac{(\theta_{k\zeta} - \widehat{\theta}_{k\zeta})^2}{\rho_{k\zeta}} \nonumber
\end{align}
where ${\mathbf A} =
{diag}\left(\theta_{11},\frac{\theta_{21}}{\delta_1},\frac{\theta_{31}}{\delta_1\delta_2},\cdots,
  \frac{\theta_{n1}}{\delta_1\cdots\delta_{n-1}}\right)
\in \mathbb{R}^{n \times n}$ is positive definite. Then, it follows from \eqref{EQ_nu1a_dot}, \eqref{EQ_nu2a_dot} and \eqref{EQ_nu3a_dot} that
\begin{flalign} 
			&\hspace{3.0cm}\dot{\nu}_{\rm tot} = \dot{\nu}_1 + \sum_{i = 2}^{n - 1} \dot{\nu}_i + \dot{\nu}_n& \nonumber
\end{flalign}
	\vspace{-0.5cm}
\begin{flalign}	
			&\hspace{3.6cm}\leqslant - \lambda_1{e}_1^2 + s_1 - \sum_{i = 2}^{n - 1}\left[\frac{\lambda_i}{\delta_1\cdots\delta_{i-1}}{e}_{i}^2
- s_{i-1} + s_{i}\right] -\frac{\lambda_n}{\delta_1\cdots\delta_{n-1}}{e}_{n}^2 - s_{n-1}& \nonumber
\end{flalign}
	\vspace{-0.5cm}
\begin{flalign}
			&\hspace{3.6cm}= - \lambda_1{e}_1^2 - \sum_{i = 2}^{n-1}\frac{\lambda_1}{\delta_1\cdots\delta_{i-1}}{e}_{i}^2 -\frac{\lambda_n}{\delta_1\cdots\delta_{n-1}}{e}_{n}^2 + \sum_{k = 1}^{n - 1} (s_k - s_k)& \nonumber
\end{flalign}
	\vspace{-0.5cm}
\begin{flalign}
			&\hspace{3.6cm}= -\pmb{e}^T{\mathbf B}\pmb{e}& \nonumber
\end{flalign} 
where ${\mathbf B} = {diag}\left(\lambda_{1},\frac{\lambda_{2}}{\delta_1},\frac{\lambda_{3}}{\delta_1\delta_2},\cdots,\frac{\lambda_{n}}{\delta_1\cdots\delta_{n-1}}\right) \in \mathbb{R}^{n \times n}$ is positive definite and every stability connector $s_k$~is~canceled by its negative counterpart $-s_k$, $\forall k \in \{1,2,...,n-1\}$. By \cite[Thm. 8.4]{Khalil2002} both the control errors and the parameter estimation errors are bounded, and $\pmb{e}(t)^T{\mathbf B}\pmb{e}(t) \to 0$ globally as $t \to \infty$, which by the positive-definiteness of $\mathbf{B}$ is equivalent~to $\pmb{e}(t)\to \bf 0$ as $t\to\infty$, i.e., $e_k(t) \to 0$, $\forall k \in \{1,2,\ldots,n\}$ as $t\to\infty$.
\end{IEEEproof}
\end{theorem}

Finally, motivated by the original concept of \textit{virtual stability} \cite[Sect. 2.9]{Zhu_VDC}, Definition \ref{def:Vstab} generalizes the results in Lemmas \ref{lem:ss_1a}--\ref{lem:ss_na} for \textit{virtual stability} of the \textit{k}th subsystem. 
\begin{defin} \label{def:Vstab}
The \textit{k}th subsystem, $\forall k \in \{1,...,n\}$, in \eqref{EQ_dyn1}--\eqref{EQ_dyn3}, combined with its respective control in \eqref{EQ_ctrl1a}--\eqref{PA_thetak}, is said to be \textit{virtually stable} if the derivative of a quadratic function
$\nu_k = \alpha_k{e}_k^2 + (\pmb{\theta}_k - \widehat{\pmb{\theta}}_k)^T{\pmb \Gamma}_k(\pmb{\theta}_k - \widehat{\pmb{\theta}}_k)$ along the trajectories~of~the error dynamics satisfies $\dot{\nu}_k \leqslant - \beta_k{e}_k^2 - s_{k-1} + s_k$ for some $\alpha_k,\beta_k > 0$ and positive-definite ${\pmb \Gamma}_k \in \mathbb{R}^{k\times k}$, where $s_{k-1}$ and $s_k$ are the stability connectors by Def.~\ref{def:stab_con} such that $s_0 = 0$ and $s_n = 0$.
\end{defin}

\begin{remark} \label{remark:virtual_stability}
Definition \ref{def:Vstab} provides generic tools to design local subsystem-based control for 
SFF systems. As we demonstrated in Theorem \ref{thm:stability}, 
\textit{virtual stability} of every SS in the sense of Definition \ref{def:Vstab} (derived from Lemmas 
\ref{lem:ss_1a}--\ref{lem:ss_na}) guarantees \textit{global asymptotic stability} of the overall system.
\end{remark}

\section{Numerical Validation}
\label{sec:validation}

In order to validate the proposed method, we consider the 3rd order nonlinear system from \cite{Yip1998,Wang2005}, namely,
\begin{empheq}[left=\empheqlbrace]{align}
  \dot{x}_1 &= a_1x_1^3 + x_2 \nonumber\\%\label{EQ_BS1}\\
  \dot{x}_2 &= a_2(x_1^2 + x_2^2) + x_3 	\label{EQ_BS2}\\
  \dot{x}_3 &= u. \nonumber%\label{EQ_BS3} 
\end{empheq}

Using \eqref{EQ_ctrl1a}--\eqref{EQ_ctrl3a}, the proposed adaptive SBC for the system in \eqref{EQ_BS2} can be designed as
\begin{empheq}[left=\empheqlbrace]{align}
  x_{\rm 2d} &= {\mathbf Y}_{1}\widehat{\pmb \theta}_1 + \lambda_1{e}_1	\nonumber\\%\label{EQ_ctrlBS1}\\
	x_{\rm 3d} &= {\mathbf Y}_{2}\widehat{\pmb \theta}_2 + \lambda_2{e}_2 + \delta_1{e}_1 \label{EQ_ctrlBS2}\\
  u &= {\rm Y}_{3}\widehat{\theta}_3 + \lambda_3{e}_3 + \delta_2{e}_2 \nonumber%\label{EQ_ctrlBS3} 
\end{empheq}
where ${\mathbf Y}_{1} = [\dot{x}_{\rm 1d}\ \ -x_1^3]$, ${\pmb \theta}_1 = [1\ \ \theta_{12}]^T$, 
${\mathbf Y}_{2} = [\dot{x}_{\rm 2d}$ $-(x_1^2 + x_2^2)]$, ${\pmb \theta}_2 = [1\ \ 
\theta_{22}]^T$, ${\rm Y}_{3} = \dot{x}_{\rm 3d}$ and ${\theta}_3 =$ ${\theta}_{31} = 1$. The 
parameters in $\widehat{\pmb \theta}_k$, $\forall k \in \{1,2,3\}$, are updated with $\mathcal{P}_k$ 
in Definition \ref{def:PA3}. For simplicity, only parameters $\theta_{12}$ and $\theta_{22}$ 
(corresponding ${a}_1$ and ${a}_2$ in the plant) are adapted in the experiments, although 
possibility to adapt $\theta_{k1}$ = 1, $\forall k \in \{1,2,3\}$, remains.

To study the global asymptotic convergence suggested by Theorem \ref{thm:stability}, the following piecewise differentiable and sufficiently smooth reference trajectory $x_{\rm 1d}(t)$ is used
\begin{align} \nonumber
x_{\rm 1d}(t) = \left\{
  \begin{array}{l}
    \hspace{-0.2cm}{\rm sin}(2\pi t){\rm tanh}(t^3),\hspace{0.2cm} {\rm if}\ 0 \leqslant t \leqslant 5\\
		\hspace{-0.2cm}{\rm sin}(2\pi t){\rm tanh}(t^3)[1-{\rm tanh}((t-5)^3)],\hspace{0.2cm} {\rm if}\ t > 5.
		\end{array} \right.
\end{align}

\begin{figure}[t]
	\vspace{0.0cm}
	\centering
	\includegraphics[width=0.65\columnwidth]{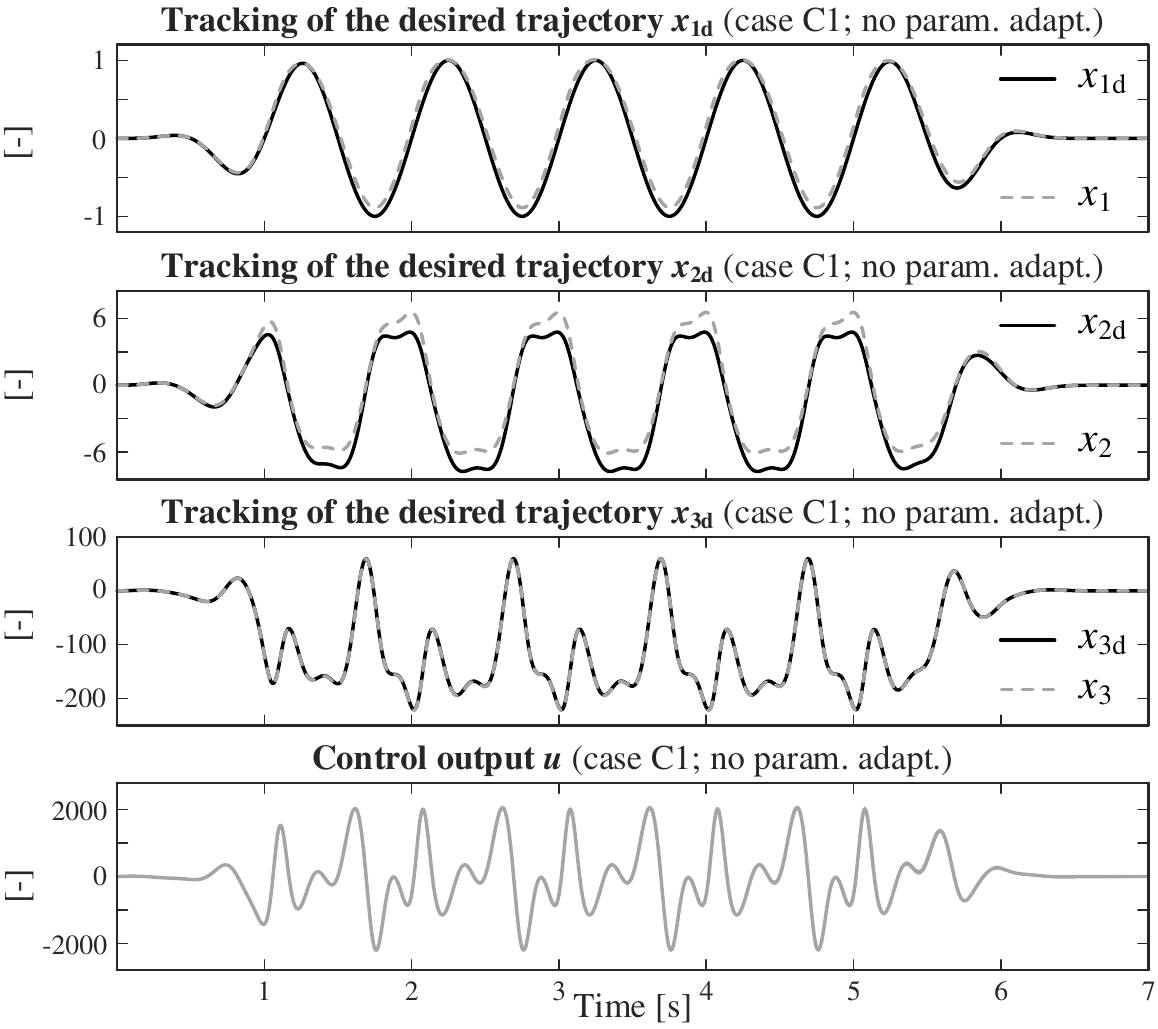}
	\vspace{-0.3cm}
	\caption{Control performance in C1 with inaccurate parameter values $\theta_{12}$ = 6 and $\theta_{22}$~=~4 in relation to the actual plant parameters $a_1$ = 5 and $a_2$ = 5. The desired trajectories are shown in black and their controlled variables in gray (plots 1--3). The last plot shows the control output $u$.} 
	\label{fig:noPA_tracking}
	\vspace{-0.0cm}
\end{figure}

\begin{figure}[t]
	\vspace{-0.1cm}
	\centering
	\includegraphics[width=0.65\columnwidth]{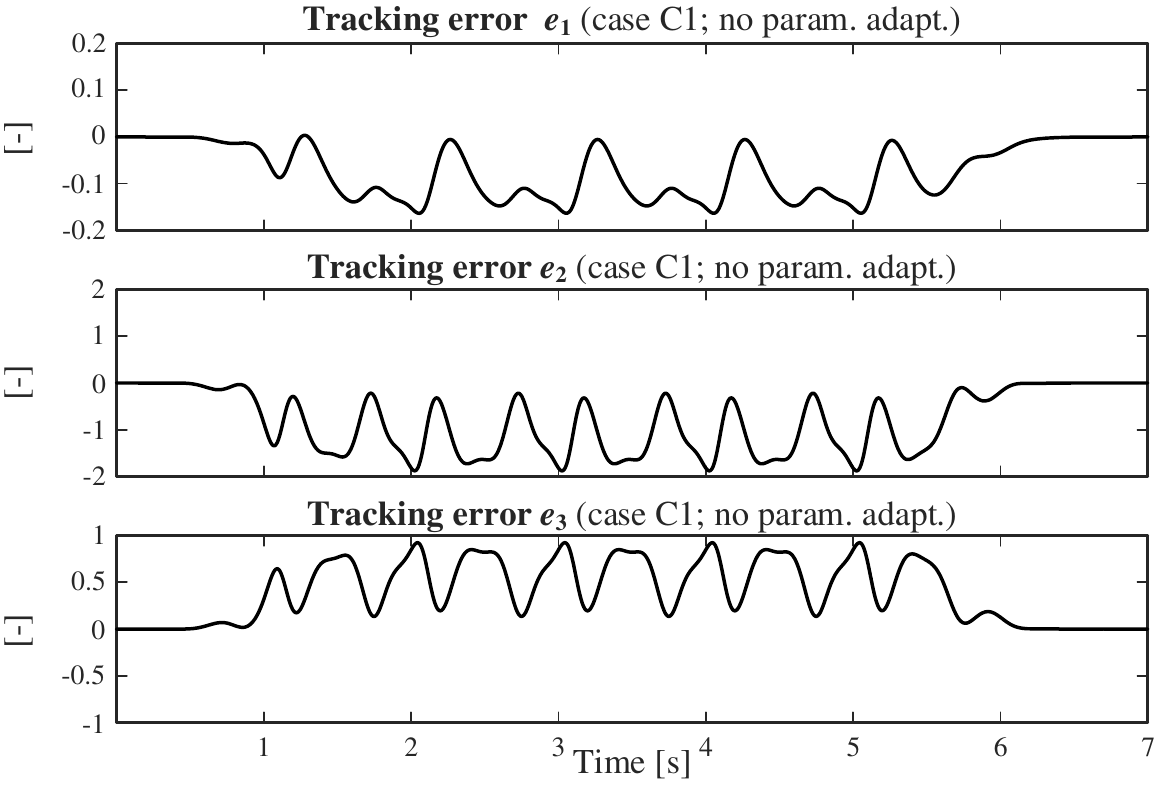}
	\vspace{-0.3cm}
	\caption{Tracking errors $e_k$, $\forall k \in \{1,2,3\}$, in C1 (adaptive control disabled).} 
	\label{fig:noPA_tracking_errors}
	\vspace{-0.0cm}
\end{figure}

\begin{figure}[t]
	\vspace{0.0cm}
	\centering
	\includegraphics[width=0.65\columnwidth]{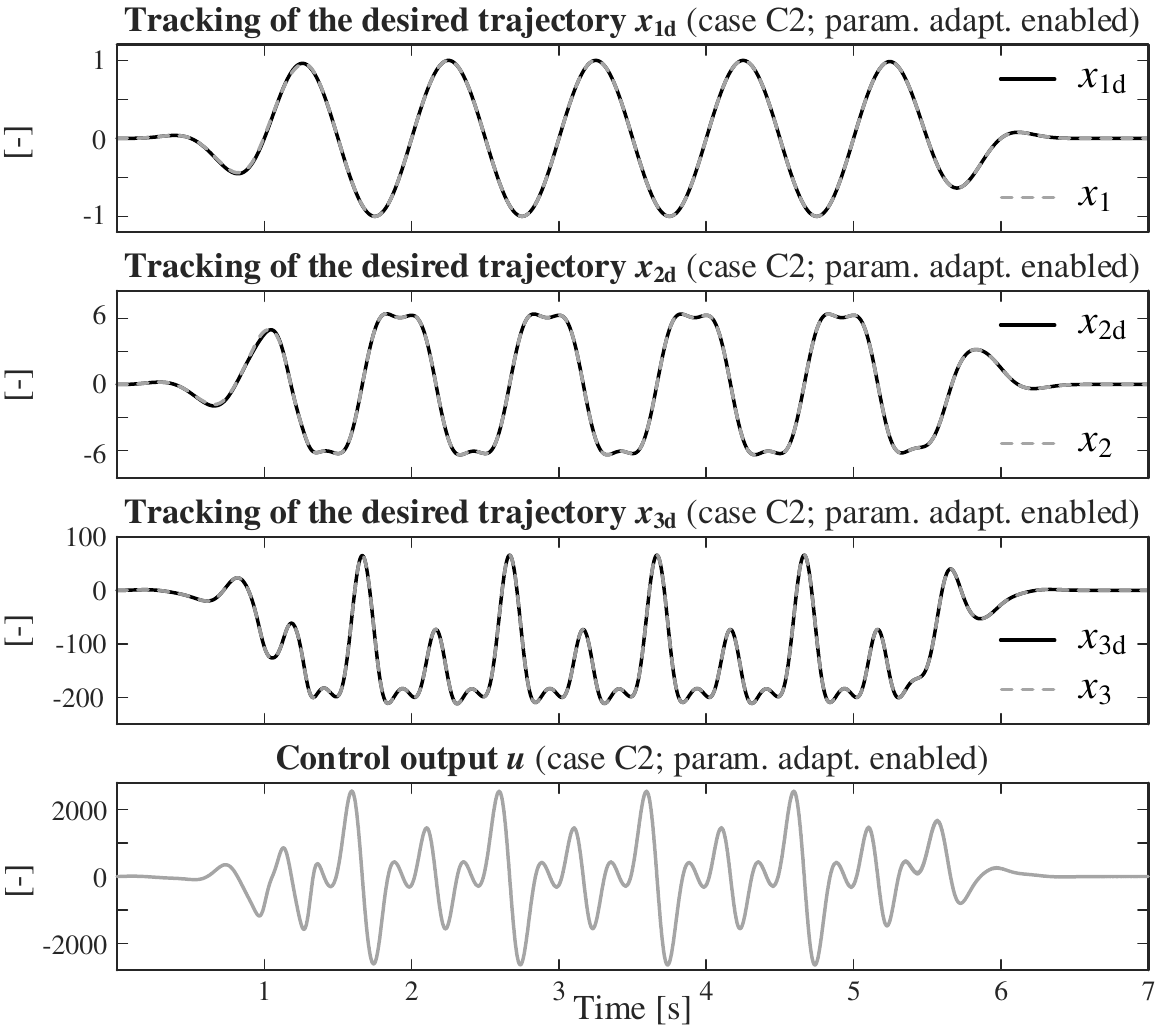}
	\vspace{-0.3cm}
	\caption{Control performance in C2 with initial~parameter values $\widehat{\theta}_{12}(0)$~=~6 and $\widehat{\theta}_{22}(0)$ = 4, while $a_1$ = 5 and $a_2$ = 5 hold for the respective plant parameters. The desired trajectories are shown in black and their~controlled variables in gray (plots 1--3). The last plot shows the control output~$u$.} 
	\label{fig:adaptive_tracking}
	\vspace{-0.0cm}
\end{figure}

\begin{figure}[t]
	\vspace{-0.1cm}
	\centering
	\includegraphics[width=0.65\columnwidth]{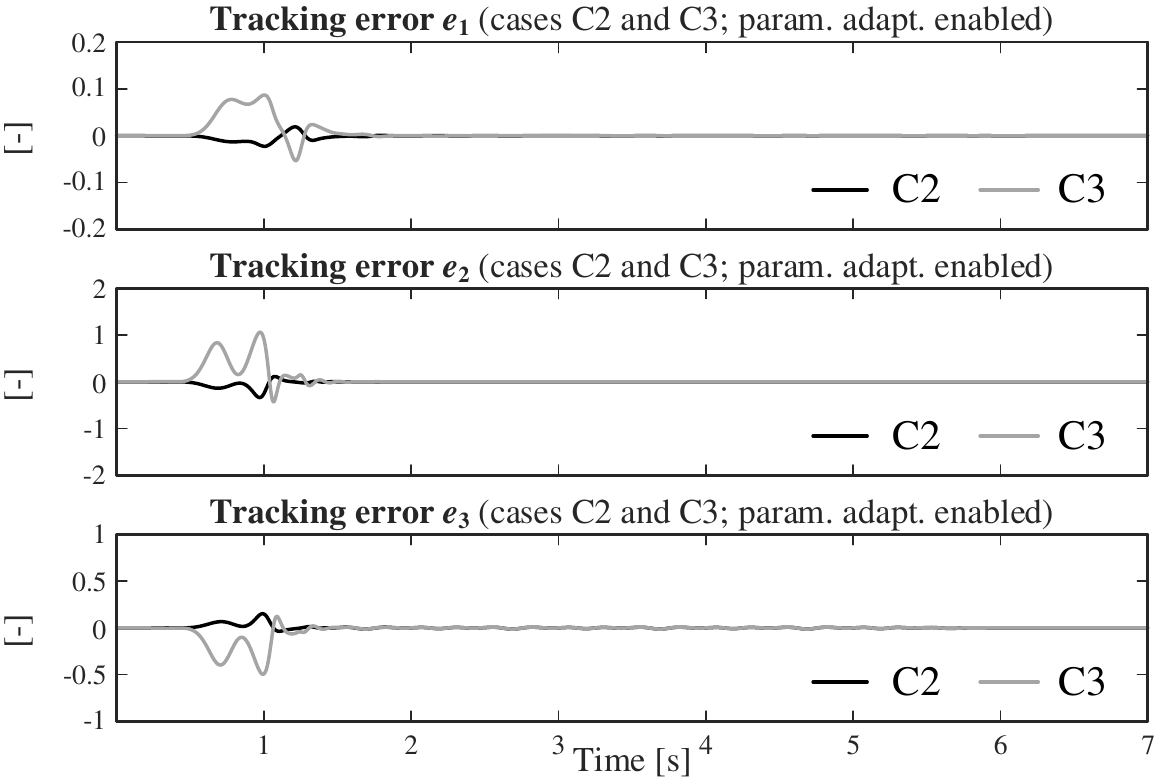}
	\vspace{-0.3cm}
	\caption{Tracking errors $e_k$, $\forall k \in \{1,2,3\}$, in C2 and C3 (adaptive~control~enabled). The results in C2 are~in~black and the results in C3 are in gray.} 
	\label{fig:adaptive_tracking_errors}
	\vspace{-0.0cm}
\end{figure}

\begin{figure}[t]
	\vspace{0.0cm}
	\centering
	\includegraphics[width=0.65\columnwidth]{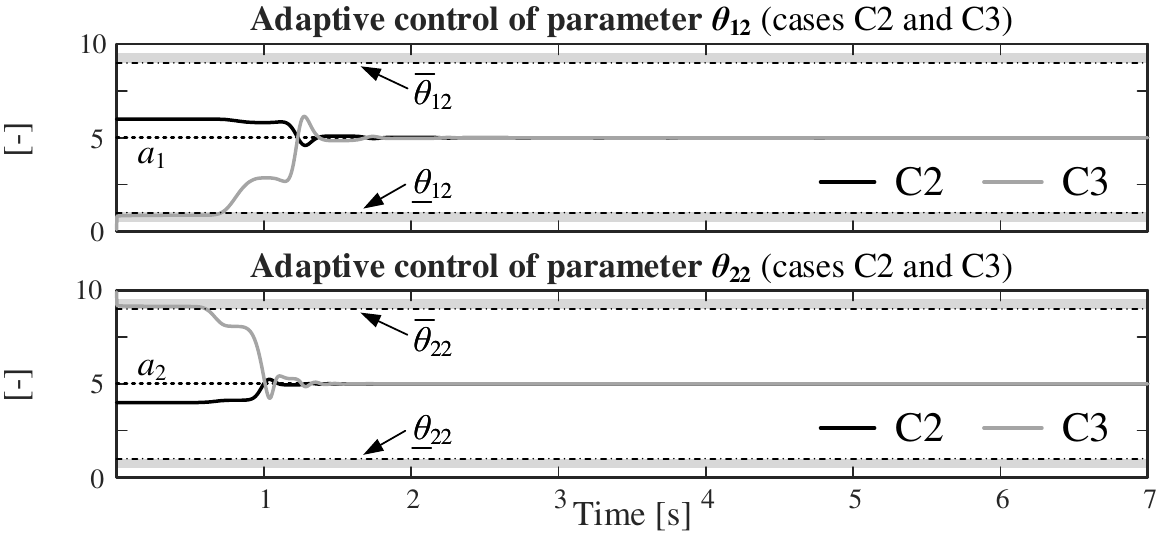}
	\vspace{-0.3cm}
	\caption{Adapted parameters $\widehat{\theta}_{12}(t)$ (the 1st plot) and $\widehat{\theta}_{22}(t)$ (the 2nd plot).~The results in C2 ($\widehat{\theta}_{12}(0)$ = 6 and $\widehat{\theta}_{22}(0)$ = 4) are given in black, while the result in C3 ($\widehat{\theta}_{12}(0)$ = 0.1 and $\widehat{\theta}_{22}(0)$ = 9.9) are given in dark gray. The plant~parameters $a_1$ = 5 and $a_2$ = 5 are shown in dashed line. The upper bound $\overline{\theta}$ = 9 and lower bound $\underline{\theta}$ = 1 are shown in dashed-dot line. The bound activation intervals (defined by $c_{\rm 21} = 0.5$ and $c_{\rm 22} = 0.5$) are in gray.} 
	\label{fig:adapted_parameters1}
	\vspace{-0.1cm}
\end{figure}

Throughout the simulations, $a_1$ = 5 and $a_2$ = 5 are used for the plant in \eqref{EQ_BS2}, and the FB gains were loosely tuned to $\lambda_1$ = 10, $\lambda_2$ = 20, $\lambda_3$ = 40, $\delta_1$ = 10 and $\delta_2$ = 20. The sample time in simulations was set to 0.01 ms to address the exponential rate of dynamics. The following three test cases are studied:
\begin{description}
\item [C1:] The baseline SBC (in Sec. \ref{subsec:SBC_design}) is employed, i.e., ${\pmb \theta}_1$, ${\pmb \theta}_2$ and ${\theta}_3$ (instead of $\widehat{\pmb \theta}_1$, $\widehat{\pmb \theta}_2$ and $\widehat{\theta}_3$) are used in \eqref{EQ_ctrlBS2}. In addition, inaccurate~FF~parameters $\theta_{12}$ = 6 and $\theta_{22}$~=~4~are used in relation to their respective plant parameters $a_1 = 5$ and $a_2 = 5$ in \eqref{EQ_BS2}. Figs. \ref{fig:noPA_tracking} and \ref{fig:noPA_tracking_errors} show the results.
\item [C2:] The proposed adaptive SBC (in Sec. \ref{subsec:ASBC_design}) is employed with initial parameter estimate values $\widehat{\theta}_{12}(0)$ = 6~and $\widehat{\theta}_{22}$(0) = 4 in \eqref{EQ_ctrlBS2}. Figs. \ref{fig:adaptive_tracking}--\ref{fig:adapted_parameters1} show the results.
\item [C3:] The proposed adaptive SBC (in Sec. \ref{subsec:ASBC_design}) is employed with initial parameter estimate values $\widehat{\theta}_{12}(0)$ = 0.1 and $\widehat{\theta}_{22}$(0) = 9.9. Figs. \ref{fig:adaptive_tracking_errors} and \ref{fig:adapted_parameters1} show the results.
\end{description}
In cases C2 and C3 with the adaptive control, the parameter update gains were set to $\rho_{12} = 
1000$, $\sigma_{12} = 1000/\rho_{12}$, $\rho_{22} = 2$ and $\sigma_{22} = 1000/\rho_{22}$; the 
parameter bounds were set to $\overline{\theta}_{12} = 9$, $\underline{\theta}_{12} = 1$, 
$\overline{\theta}_{22} = 9$ and $\underline{\theta}_{22} = 1$; and the activation intervals beyond 
the bounds were set to 0.5 ($c_{\rm 21} = 0.5$ and $c_{\rm 22} = 0.5$)

Figs. \ref{fig:noPA_tracking} and \ref{fig:noPA_tracking_errors} show the results in case C1,~where~the~control for the plant is designed by the theory in Section \ref{subsec:SBC_design}. However, inaccurate FF parameters (i.e., $\theta_{12} \neq a_1$ and $\theta_{22} \neq a_2$)~are used such that, in cases C2 and C3, comparisons can be made to the proposed adaptive SBC. In plots 1--3, Fig. \ref{fig:noPA_tracking} shows the desired trajectory $x_{k{\rm d}}$, $\forall k \in \{1,2,3\}$, in black and its controlled state $x_{k}$ in gray. The last plot shows the control output $u$. The detailed tracking errors are shown in Fig. \ref{fig:noPA_tracking_errors}, the maximum absolute tracking errors being $|e_1|_{\rm max}$ = 0.163, $|e_2|_{\rm max}$ = 1.877 and $|e_3|_{\rm max}$ = 0.924. As can be seen, noticeable tracking errors occur in the transition phases due to the parametric uncertainty.

\textit{Figs. \ref{fig:adaptive_tracking}--\ref{fig:adapted_parameters1} show the main results of the study}~with~the~proposed adaptive SBC in \eqref{EQ_ctrlBS2}. Fig. \ref{fig:adaptive_tracking} shows the tracking results in case C2 where the initial values for the parameter estimates are selected in accordance to case C1, i.e., $\widehat{\theta}_{12}(0)$ = 6 and $\widehat{\theta}_{22}(0)$ = 4. As the black lines in Fig. \ref{fig:adaptive_tracking_errors} shows, the tracking errors~are substantially decreased in relation to case C1, with the maximum absolute tracking errors $|e_1|_{\rm max}$ = 0.023, $|e_2|_{\rm max}$ = 0.336 and $|e_3|_{\rm max}$ = 0.152. As predicted by the theory, global asymptotic convergence is achieved. Fig. \ref{fig:adapted_parameters1} shows the behavior of the parameter estimates $\widehat{\theta}_{12}$ and $\widehat{\theta}_{22}$ in black, illustrating that the proposed projection function ${\mathcal P}_k$ actively pushes the parameter values toward their real values in the plant.

In the last case C3, the initial parameter values are~set~outside the projection function ${\mathcal P}_k$ bounds such that $\widehat{\theta}_{12}$(0) = 0.1 and $\widehat{\theta}_{22}$(0) = 9.9. The results are shown~in~Figs.~\ref{fig:adaptive_tracking_errors}~and~\ref{fig:adapted_parameters1} in gray. Despite a significant inaccuracy in the initial parameter values, the projection function ${\mathcal P}_k$ actively pushes the parameter values toward their real values in the plant (see Fig. \ref{fig:adapted_parameters1}), with the maximum absolute tracking errors $|e_1|_{\rm max}$ = 0.087, $|e_2|_{\rm max}$ = 1.060 and $|e_3|_{\rm max}$ = 0.496 (see Fig. \ref{fig:adaptive_tracking_errors}). After 1.5 s the control behavior in case C3 becomes virtually identical to case C2.

\section{Conclusions}
\label{sec:conclusions}

This study proposed an \textit{adaptive subsystem-based control} for controlling \textit{n}th-order 
SFF systems with parametric uncertainties. As an alternative for backstepping, we provided 
systematic~and straightforward tools for globally~asymptotically stable control while avoiding a 
growth of the control design complexity when the system order \textit{n} increases. The proposed 
method is modular in the sense that the control for every SS can be designed with a single 
generic-form equation such that changing SS dynamics or removing/adding SSs do not affect to the 
control laws in the remaining SSs. For the method, we reformulated the original concept of \textit{virtual 
stability} in \cite[Def. 2.17]{Zhu_VDC} and proposed a specific \textit{stability connector} to address dynamic 
interactions between the adjacent SSs. These features enable that both the control design and the stability 
analysis can be performed locally at a SS level (as opposed to the whole system); see 
Remark \ref{remark:virtual_stability}. We proposed also a smooth projection function ${\mathcal P}_k$ 
for the system parametric uncertainties. Theoretical developments on global asymptotic 
convergence (in Theorem \ref{thm:stability}) were verified in numerical simulations. Semi-SFF 
systems with unknown dynamics remain a subject for future studies. 

\appendices  

\section{The Projection Function $\mathcal{P}_{k}$}
\label{App:projF}

Consider the piecewise-continuous projection function $\mathcal{P}_{k}$ in Definition \ref{def:PA3}. 
Parameters $a$ and $b$ define the lower and upper bounds for $\mathcal{P}_{k}$ such that $b 
\geqslant a > 0$.~Within~the bounds, $\dot{\mathcal{P}}_{k}$ is driven by $\rho {\rm p}(t)$ and the behavior of $\mathcal{P}_{k}$ is equal to $\mathcal{P} \in C^1$ and $\mathcal{P}_2 \in C^2$ in \cite{Zhu_VDC,Zhu1999adaptive}. Outside 
the bounds, a corrective term $\sigma{\kappa}$ is designed to bring $\mathcal{P}_{k}$ back 
toward~the bounds. The parameter $c$ defines activation interval lengths  $(a-c,a)$ and $(b,b+c)$ 
for the switching functions $S_a(\mathcal{P}_k)$ and $S_b(\mathcal{P}_k)$. 

Let $|{\rm p}^{(n-k)}(t)| < +\infty$, $\forall k \in \{1,...,n\}$. To guarantee the existence of 
${\mathcal{P}^{(k)}_{k}}$ outside the bounds, the switching functions $S_a(\mathcal{P}_{k}): (a-c, 
a) \to (1, 0)$ and $S_b(\mathcal{P}_{k}): (b, b+c) \to (0, 1)$ are required to satisfy the boundary 
conditions [in \eqref{EQ_BCa} and \eqref{EQ_BCb}]
\begin{equation}
\begin{aligned} \label{EQ_BCa}
    &{\lim_{x \to (a-c)^+}}S_a(x) = 1,\hspace{0.2cm} {\lim_{x \to (a-c)^+}}S_a^{(j)}(x) = 0, \\[-0.15cm]
		&\hspace{0.25cm}{\lim_{x \to a^-}}S_a(x) = 0\ \ {\rm and}\ {\lim_{x \to a^-}}S_a^{(j)}(x) = 0, 
\end{aligned} 
\end{equation}	
\begin{equation}
\begin{aligned} \label{EQ_BCb}
    &{\lim_{x \to (b+c)^-}}S_b(x) = 1,\hspace{0.2cm} {\lim_{x \to (b+c)^-}}S_b^{(j)}(x) = 0, \\[-0.15cm]
		&\hspace{0.25cm}{\lim_{x \to b^+}}S_b(x) = 0\ \ {\rm and}\ {\lim_{x \to b^+}}S_b^{(j)}(x) = 0,
\end{aligned} 
\end{equation}	
$\forall j \in \{1,...,n-1\}$. Definition \ref{Switching_functions} provides smooth and strictly decreasing solution for $S_a(\mathcal{P}_{k})$, satisfying \eqref{EQ_BCa}, and smooth and strictly increasing solution for $S_b(\mathcal{P}_{k})$, satisfying \eqref{EQ_BCb}.

\begin{defin} \label{Switching_functions}
$S_a(\mathcal{P}_{k}): (a-c, a) \to (1, 0)$ is a \textit{smooth and strictly decreasing} switching function defined as
\begin{flalign} \nonumber%\label{EQ_Sa}
  S_a(\mathcal{P}_{k}) := \frac{1}{2}\left[ 1-\tanh\left( \frac{1}{a-c-\mathcal{P}_{k}} + 
  \frac{1}{a-\mathcal{P}_{k}} \right) \right]
\end{flalign}
and $S_b(\mathcal{P}_{k}): (b, b+c) \to (0, 1)$ is a \textit{smooth and strictly increasing} switching function defined as
\begin{flalign} \nonumber%\label{EQ_Sb}
  S_b(\mathcal{P}_{k}) := \frac{1}{2}\left[ 1+\tanh\left( \frac{1}{b-\mathcal{P}_{k}} + 
  \frac{1}{b+c-\mathcal{P}_{k}} \right) \right].
\end{flalign}
\end{defin}

The projection function in~\eqref{EQ_PAdef3} has the following property.
\begin{lemma} \label{lem:ProjF3}
For any constant $\mathcal{P}_{c}$ with $a \leqslant \mathcal{P}_{c} \leqslant b$ we have
\begin{equation} \label{EQ_PA_ineq3}
({\mathcal{P}_{c}} - \mathcal{P}_{k})\left({\rm p}(t) - \frac{1}{\rho}\dot{\mathcal{P}}_{k}\right) \leqslant -\sigma\kappa^2 \leqslant 0.
\end{equation}
\end{lemma}

\begin{IEEEproof}
The proof of Lemma \ref{lem:ProjF3} follows a similar procedure as the proof of Lemma 2.10 in 
\cite{Zhu_VDC}. Let $a \leqslant {\mathcal{P}_{c}} \leqslant b$. Then, for a constant 
${\mathcal{P}_{c}}$ we have 
\begin{align}
	({\mathcal{P}_{c}} - \mathcal{P}_{k}) &\geqslant (a - \mathcal{P}_{k}) \label{EQ_projF_ineq_a1}\\
	(\mathcal{P}_{k} - {\mathcal{P}_{c}}) &\geqslant (\mathcal{P}_{k} - b). \label{EQ_projF_ineq_b1}
\end{align}

Substituting \eqref{EQ_PAdef3} into \eqref{EQ_PA_ineq3} we get
\begin{equation} \label{EQ_ineg_proof1}
({\mathcal{P}_{c}} - \mathcal{P}_{k})\left({\rm p}(t) - \frac{1}{\rho}\dot{\mathcal{P}}_{k}\right) = -\sigma({\mathcal{P}_{c}} - \mathcal{P}_{k}){\kappa}. 
\end{equation}

When $\mathcal{P}_{n} \leqslant a - c$, Definition~\ref{def:PA3} yields ${\kappa} = (a - 
\mathcal{P}_{k})$. Using \eqref{EQ_projF_ineq_a1},~we~get
\begin{align} \label{EQ_ineg_proof2}
	-\sigma({\mathcal{P}_{c}} - \mathcal{P}_{k}){\kappa} \leqslant -\sigma(a - \mathcal{P}_{k}){\kappa} =-\sigma\kappa^2 \leqslant 0.
\end{align}

When $a - c < \mathcal{P}_{k} < a$, Definition~\ref{def:PA3} yields ${\kappa} = (a - 
\mathcal{P}_{k})S_a(\mathcal{P}_{k})$ and $S_a(\mathcal{P}_{k}) \in (0,1)$. Using 
\eqref{EQ_projF_ineq_a1},~we~get 
\begin{align} \label{EQ_ineg_proof3}
	-\sigma({\mathcal{P}_{c}} - \mathcal{P}_{k}){\kappa} &\leqslant -\sigma(a - \mathcal{P}_{k}){\kappa} \nonumber\\
					&<-\sigma(a - \mathcal{P}_{k})S_a(\mathcal{P}_{k}){\kappa} \nonumber\\
					&=-\sigma\kappa^2 \leqslant 0.
\end{align}

When $a \leqslant \mathcal{P}_{k} \leqslant b$, Definition~\ref{def:PA3} yields ${\kappa} = 0$ and we get
\begin{align} \label{EQ_ineg_proof4}
	-\sigma({\mathcal{P}_{c}} - \mathcal{P}_{k}){\kappa} = 0.
\end{align}

When $b < \mathcal{P}_{k} < b + c$, Definition~\ref{def:PA3} yields ${\kappa} = (b - 
\mathcal{P}_{k})S_b(\mathcal{P}_{k})$ and $S_b(\mathcal{P}_{k}) \in (0,1)$. Using 
\eqref{EQ_projF_ineq_b1},~we~get
\begin{align} \label{EQ_ineg_proof5}
	-\sigma({\mathcal{P}_{c}} - \mathcal{P}_{k}){\kappa} &= \sigma(\mathcal{P}_{k} - {\mathcal{P}_{c}}){\kappa} \nonumber\\
					&\leqslant \sigma(\mathcal{P}_{k} - b){\kappa} \nonumber\\
					&< \sigma(\mathcal{P}_{k} - b)S_b(\mathcal{P}_{k}){\kappa} \nonumber\\
					&= - \sigma(b - \mathcal{P}_{k})S_b(\mathcal{P}_{k}){\kappa} \nonumber\\
					&= - \sigma\kappa^2\leqslant 0.
\end{align}

When $\mathcal{P}_{k} \geqslant b + c$, Definition~\ref{def:PA3} yields ${\kappa} = (b - 
\mathcal{P}_{k})$. Using \eqref{EQ_projF_ineq_b1}, we get
\begin{align} \label{EQ_ineg_proof6}
	-\sigma({\mathcal{P}_{c}} - \mathcal{P}_{k}){\kappa} &= \sigma(\mathcal{P}_{k} - {\mathcal{P}_{c}}){\kappa} \nonumber\\
					&\leqslant \sigma(\mathcal{P}_{k} - b){\kappa} \nonumber\\
					&= - \sigma(b - \mathcal{P}_{k}){\kappa} \nonumber\\
					&= - \sigma\kappa^2\leqslant 0.
\end{align}
Finally, \eqref{EQ_ineg_proof2}--\eqref{EQ_ineg_proof6} together with \eqref{EQ_ineg_proof1} 
complete the proof. 
\end{IEEEproof}

\section{Proofs of Lemmas \ref{lem:ss_1a}, \ref{lem:ss_ia}, and~\ref{lem:ss_na}}
\label{App:lemmas_ss}

\textit{Proof of Lemma~\textup{\ref{lem:ss_1a}}}: Using \eqref{EQ_error_dyn1a}, Definition \ref{def:stab_con}~and~Lemma \ref{lem:ProjF3}, the derivative of the quadratic function $\nu_1$ in \eqref{EQ_nu1a} can be written as
\begin{align}
  \dot{\nu}_{1} &= {e}_1\theta_{11}\dot{e}_1 - \sum_{\zeta = 1}^{j}(\theta_{1\zeta} - \widehat{\theta}_{1\zeta})\frac{\dot{\widehat{\theta}}_{1\zeta}}{\rho_{1\zeta}}\nonumber\\
  &= - \lambda_1{e}_1^2 + g_1({x}_1){e}_1{e}_2 + {e}_1{\mathbf Y}_{1}({\pmb \theta}_1 - \widehat{\pmb \theta}_1) - \sum_{\zeta = 1}^{j}(\theta_{1\zeta} - \widehat{\theta}_{1\zeta})\frac{\dot{\widehat{\theta}}_{1\zeta}}{\rho_{1\zeta}} \nonumber\\
  &= - \lambda_1{e}_1^2 + g_1({x}_1){e}_1{e}_2 + \sum_{\zeta=1}^{j}({\theta}_{1\zeta} - \widehat{{\theta}}_{1\zeta})\left({\rm p}_{1\zeta} - \frac{\dot{\widehat{{\theta}}}_{1\zeta}}{{\rho}_{1\zeta}}\right) \nonumber\\
  &\leqslant - \lambda_1{e}_1^2 + s_1 \nonumber
\end{align} which completes the proof of Lemma \ref{lem:ss_1a}.  \hfill $\blacksquare$

\textit{Proof of Lemma \textup{\ref{lem:ss_ia}}}: Using \eqref{EQ_error_dyn2a}, Definition \ref{def:stab_con}~and~Lemma \ref{lem:ProjF3}, the derivative of the quadratic function $\nu_i$ in \eqref{EQ_nu2a}, $\forall i \in \{2,...,$ $n-1\}$, can be written as
\begin{flalign} 
		&\hspace{3.3cm}\dot{\nu}_i = {e}_i\frac{\theta_{i1}}{\delta_1\cdots\delta_{i-1}}\dot{e}_{i} - \frac{1}{\delta_1\cdots\delta_{i-1}}\sum_{\zeta = 1}^{j}(\theta_{i\zeta} - \widehat{\theta}_{i\zeta})\frac{\dot{\widehat{\theta}}_{i\zeta}}{\rho_{i\zeta}}& \nonumber
\end{flalign}
	\vspace{-0.5cm}
\begin{flalign} 
		&\hspace{3.6cm}= {e}_{i}\frac{1}{\delta_1\cdots\delta_{i-1}}\Big[g_i(\pmb{x}_i)e_{i+1} - \delta_{i-1}g_{i-1}(\pmb{x}_{i-1}){e}_{i-1} - \lambda_i{e}_i + {\mathbf Y}_i({\pmb \theta}_i - \widehat{\pmb \theta}_i)\Big]& \nonumber
\end{flalign}
	\vspace{-0.5cm}
\begin{flalign} 
		&\hspace{4.0cm} - \frac{1}{\delta_1\cdots\delta_{i-1}}\sum_{\zeta = 1}^{j}(\theta_{i\zeta} - \widehat{\theta}_{i\zeta})\frac{\dot{\widehat{\theta}}_{i\zeta}}{\rho_{i\zeta}}& \nonumber
\end{flalign}
	\vspace{-0.5cm}
\begin{flalign} 
		&\hspace{3.6cm}= - \frac{\lambda_i}{\delta_1\cdots\delta_{i-1}}{e}_{i}^2 - \frac{1}{\delta_1\cdots\delta_{i-2}}g_{i-1}(\pmb{x}_{i-1}){e}_{i-1}{e}_i + \frac{1}{\delta_1\cdots\delta_{i-1}}g_i(\pmb{x}_i){e}_i{e}_{i+1}& \nonumber
\end{flalign}
	\vspace{-0.5cm}
\begin{flalign}  
		&\hspace{4.0cm}+ \frac{1}{\delta_1\cdots\delta_{i-1}}\sum_{\zeta=1}^{j}({\theta}_{i\zeta} - \widehat{{\theta}}_{i\zeta})\left({\rm p}_{i\zeta} - \frac{\dot{\widehat{{\theta}}}_{i\zeta}}{{\rho}_{i\zeta}}\right)& \nonumber
\end{flalign}
	\vspace{-0.5cm}
\begin{flalign} 
		&\hspace{3.6cm}\leqslant - \frac{\lambda_i}{\delta_1\cdots\delta_{i-1}}{e}_i^2 - s_{i-1} + s_{i}& \nonumber
\end{flalign}
which completes the proof of Lemma \ref{lem:ss_ia}.  \hfill
$\blacksquare$

\textit{Proof of Lemma~\textup{\ref{lem:ss_na}}}: Using \eqref{EQ_error_dyn3a}, Definition \ref{def:stab_con}~and~Lemma \ref{lem:ProjF3}, the derivative of the quadratic function $\nu_n$ in \eqref{EQ_nu3a} can be written as
\begin{align}
  \dot{\nu}_n &= {e}_{n}\frac{\theta_{n1}}{\delta_1\cdots\delta_{n-1}}\dot{e}_{n} - \frac{1}{\delta_1\cdots\delta_{n-1}}\sum_{\zeta = 1}^{j}(\theta_{n\zeta} - \widehat{\theta}_{n\zeta})\frac{\dot{\widehat{\theta}}_{n\zeta}}{\rho_{n\zeta}} \nonumber\\ 
							&= {e}_{n}\frac{1}{\delta_1\cdots\delta_{n-1}}\Big[-\lambda_n{e}_{n} - \delta_{n-1}g_{n-1}(\pmb{x}_{n-1}){e}_{n-1} + {\mathbf Y}_{n}({\pmb \theta}_n - \widehat{\pmb \theta}_n)\Big] - \frac{1}{\delta_1\cdots\delta_{n-1}}\sum_{\zeta = 1}^{j}(\theta_{n\zeta} - \widehat{\theta}_{n\zeta})\frac{\dot{\widehat{\theta}}_{n\zeta}}{\rho_{n\zeta}} \nonumber\\ 
							&= -\frac{\lambda_n}{\delta_1\cdots\delta_{n-1}}{e}_{n}^2 - \frac{1}{\delta_1\cdots\delta_{n-2}}g_{n-1}(\pmb{x}_{n-1}){e}_{n-1}{e}_{n} + \frac{1}{\delta_1\cdots\delta_{n-1}}\sum_{\zeta=1}^{j}({\theta}_{n\zeta} - \widehat{{\theta}}_{n\zeta})\left({\rm p}_{n\zeta} - \frac{\dot{\widehat{{\theta}}}_{n\zeta}}{{\rho}_{n\zeta}}\right) \nonumber\\ 
							&\leqslant - \frac{\lambda_n}{\delta_1\cdots\delta_{n-1}}{e}_{n}^2 - s_{n-1} \nonumber
\end{align} 
which completes the proof of Lemma \ref{lem:ss_na}.  \hfill $\blacksquare$

\ifCLASSOPTIONcaptionsoff
  \newpage
\fi

\end{document}